# Quasi-Newton Methods: A New Direction

Philipp Hennig                                                                  philipp.hennig@tuebingen.mpg.de
Martin Kiefel                                                                   martin.kiefel@tuebingen.mpg.de
Max Planck Institute for Intelligent Systems, Spemannstraße, 72076 Tübingen, Germany

## Abstract

Four decades after their invention, quasi-Newton methods are still state of the art in unconstrained numerical optimization. Although not usually interpreted thus, these are learning algorithms that fit a local quadratic approximation to the objective function. We show that many, including the most popular, quasi-Newton methods can be interpreted as approximations of Bayesian linear regression under varying prior assumptions. This new notion elucidates some shortcomings of classical algorithms, and lights the way to a novel nonparametric quasi-Newton method, which is able to make more efficient use of available information at computational cost similar to its predecessors.

## 1. Introduction

Quasi-Newton algorithms are arguably the most popular class of nonlinear numerical optimization methods, used widely in numerical applications not just in machine learning. Their defining property is that they iteratively build estimators $B_i$ for the Hessian $B(x) = \nabla \nabla^\top f(x)$ of the objective function $f(x)$, from observations of $f$'s gradient $\nabla f(x)$, at each iteration searching for a local minimum along a line search direction $-B_i^{-1} \nabla f(x)$, an estimate of the eponymous Newton-Raphson search direction. Some of the most widely known members of this family include Broyden's (1965) method, the SR1 formula (Davidon, 1959; Broyden, 1967), the DFP method (Davidon, 1959; Fletcher & Powell, 1963) and the BFGS method (Broyden, 1969; Fletcher, 1970; Goldfarb, 1970; Shanno, 1970). Decades of continued research effort in this area make it impossible to give even a superficial overview over the available literature. The textbooks by Nocedal & Wright (1999) and Boyd & Vandenberghe (2004) are good modern starting points for readers interested in background. An insightful and extensive contemporary review was compiled by Dennis & Morée (1977). The ubiquity of optimization problems in machine learning has made these algorithms tools of the trade. But, perhaps because they predate machine learning itself, they have rarely been studied as learning algorithms in their own right. This paper offers a probabilistic analysis.

Throughout, let $f : \mathbb{R}^N \to \mathbb{R}$ be a sufficiently regular, not necessarily convex, function; $\nabla f : \mathbb{R}^N \to \mathbb{R}^N$ its gradient; $B : \mathbb{R}^N \to \mathbb{R}^{N \times N}$ its Hessian. We consider iterative algorithms moving from location $\boldsymbol{x}_{\ell-1} \in \mathbb{R}^D$ to location $\boldsymbol{x}_\ell$. The algorithm performs consecutive *line searches* along one-dimensional subspaces $\boldsymbol{x}_i(\alpha) = \alpha \boldsymbol{e}_i + \boldsymbol{x}_i^0$, with $\alpha \in \mathbb{R}_+$ and a unit length vector $\boldsymbol{e}_i \in \mathbb{R}^N$ spanning the line search space starting at $\boldsymbol{x}_i^0$. Evaluations at $\boldsymbol{x}_i$ evince the gradient $\nabla f(\boldsymbol{x}_i)$ (and usually also $f(\boldsymbol{x}_i)$, though this will not feature in this paper). The goal is to find a candidate $\boldsymbol{x}^*$ for a local minimum: a root $\nabla f(\boldsymbol{x}^*) = \boldsymbol{0}$ of the gradient.

The derivations of classical quasi-Newton algorithms proceed along the following line of argument: We require an update rule incorporating an observation $\nabla f(x_{i+1})$ into a current estimate $\hat{B}_i$ to get a new estimate $\hat{B}_{i+1}$, subject to the following desiderata:

**Low Rank/Cost Updates** Optimization problems regularly have dimensionality above $N \sim 10^3$, even beyond $N \sim 10^6$. So the update should be of low rank $M$ (usually $M = 1$ or 2), because, by Schur's lemma, it has (worst-case) cost $\mathcal{O}(N^2 + NM + M^3)$.

**Consistency with Quadratic Model** If $f$ is locally described well to second order, then

$$\boldsymbol{y}_i \equiv \nabla f(\boldsymbol{x}_i) - \nabla f(\boldsymbol{x}_{i-1}) \approx B(\boldsymbol{x}_i) \boldsymbol{s}_i, \qquad (1)$$

with $\boldsymbol{s}_i \equiv \boldsymbol{x}_i - \boldsymbol{x}_{i-1}$. Because this is the fundamental idea behind this family of algorithms, it is also known as *the quasi-Newton equation* (Dennis & Morée, 1977).

**Symmetry** The Hessian of twice differentiable functions is symmetric; so its estimator should be, too.

**Positive Definiteness** *Convex* functions have positive definite Hessians everywhere. Over time, it has become common conviction that, even for non-convex problems, positive definiteness of the estimator is desirable.

This paper's contributions are twofold: Section 2 offers a probabilistic viewpoint on classical quasi-Newton methods, in the process showing that symmetry is only achieved in a partial, definiteness in only a weak way by the classical algorithms. In Section 3 we use these insights to construct a novel nonparametric Bayesian quasi-Newton algorithm; this addresses several shortcomings of classic algorithms, and increases performance at only mildly higher cost.

We will write $\vec{X}$ to indicate an $n \times m$ matrix $X$ stacked row-wise into a vector of length $nm$. Its elements can be enumerated by an index set $(ij) \in [1,n] \times [1,m]$. The symbol $\otimes$ denotes the Kronecker product: $(a \otimes b)_{(ij)(k\ell)} = a_{ik}b_{j\ell}$. It allows a compact notation for vectorised matrices: $\overrightarrow{XYZ} = (X \otimes Z^\top)\vec{Y}$. If $A$ and $B$ have size $I \times J$ and $K \times L$, respectively, then $A \otimes B$ has size $IK \times JL$, and $\lim_{\gamma \to 0}(Y \otimes \gamma Z) = \vec{0}$ for *any* fixed, finite $Y$ and $Z$.

## 2. Quasi-Newton Methods as approximate Bayesian Regressors

From a probabilistic perspective, Equation (1) is a likelihood for $B$. Using $\boldsymbol{s}_i = \boldsymbol{x}_i - \boldsymbol{x}_{i-1}$, we can write it using Dirac's distribution

$$p(\boldsymbol{y}_i | B, \boldsymbol{s}_i) = \delta(\boldsymbol{y}_i - B\boldsymbol{s}_i) = \lim_{\beta \to 0} \mathcal{N}\left[\boldsymbol{y}_i; \mathcal{S}_i^\top \vec{B}, V_{i-1} \otimes \beta\right] \quad (2)$$

with any arbitrary $N \times N$ matrix $V_{i-1}$, a scalar $\beta$, and the linear operator $\mathcal{S}_i = (I \otimes \boldsymbol{s}_i)$. Of course, the $N$ real numbers in $\boldsymbol{y}_i$ are not sufficient to identify the $N^2$ numbers in $B$. Classical derivations (Dennis & Morée, 1977; Nocedal & Wright, 1999) thus introduce a *regularizer* based on the weighted Frobenius norm around the current best estimate $B_{i-1}$ from previous iterations. The weight in the Frobenius norm is encoded using a positive definite matrix, which we will suggestively call $V_{i-1}^{-1}$ and, without loss of generality, identify with the $V_{i-1}$ of Eq. (2)

$$\|B - B_{i-1}\|_{F,V_{i-1}^{-1}} \equiv \mathrm{tr}(V_{i-1}^{-1}(B-B_{i-1})^\top V_{i-1}^{-1}(B-B_{i-1}))$$
$$= (\vec{B} - \vec{B}_{i-1})^\top (V_{i-1}^{-1} \otimes V_{i-1}^{-1})(\vec{B} - \vec{B}_{i-1}). \quad (3)$$

The new estimate is the unique matrix $B_i$ minimizing the regularizer subject to Eq. (2). Inspecting Eq. (3) we see that, up to isomorphisms, the Frobenius regularizer is the negative logarithm of a Gaussian *prior*

$$p(B) = \mathcal{N}\left[\vec{B}; \vec{B}_{i-1}, \Sigma_{i-1} \equiv (V_{i-1} \otimes V_{i-1})\right]. \quad (4)$$

Gaussian likelihoods are conjugate to Gaussian priors. So the posterior is Gaussian, too, even for the limit case of a Dirac likelihood. A few lines of algebra[1] show that the posterior has mean and covariance

$$B_i = B_{i-1} + \frac{(\boldsymbol{y}_i - B_{i-1}\boldsymbol{s}_i)\boldsymbol{s}_i^\top V_{i-1}}{\boldsymbol{s}_i^\top V_{i-1}\boldsymbol{s}_i} \quad \text{and} \quad (5)$$

$$\Sigma_i = V_{i-1} \otimes \left(V_{i-1} - \frac{V_{i-1}\boldsymbol{s}_i \boldsymbol{s}_i^\top V_{i-1}}{\boldsymbol{s}_i^\top V_{i-1}\boldsymbol{s}_i}\right) \equiv V_{i-1} \otimes V_i, \quad (6)$$

respectively. The new mean is a rank-1 update of the old mean, and the rank of the new covariance $\Sigma_i$ is one less than that of $\Sigma_{i-1}$. The posterior mean has maximum posterior probability (minimal regularized loss), and is thus our new point estimate. Choosing a unit variance prior $\Sigma_{i-1} = \boldsymbol{I} \otimes \boldsymbol{I}$ recovers one of the oldest quasi-Newton algorithms: *Broyden's method* (1965):

$$B_i = B_{i-1} + \frac{(\boldsymbol{y}_i - B_{i-1}\boldsymbol{s}_i)\boldsymbol{s}_i^\top}{\boldsymbol{s}_i^\top \boldsymbol{s}_i} \quad (7)$$

Broyden's method does not satisfy the third requirement of Section 1: the updated estimate is, in general, not a symmetric matrix. A supposed remedy for this problem, and in fact the *only* rank-1 update rule that obeys Eq. (2) (Dennis & Morée, 1977) is the *symmetric rank 1 (SR1)* method (Davidon, 1959; Broyden, 1967):

$$B_i = B_{i-1} + \frac{(\boldsymbol{y}_i - B_{i-1}\boldsymbol{s}_i)(\boldsymbol{y}_i - B_{i-1}\boldsymbol{s}_i)^\top}{\boldsymbol{s}_i^\top (\boldsymbol{y}_i - B_{i-1}\boldsymbol{s}_i)}. \quad (8)$$

The SR1 update rule has acquired a controversial reputation (e.g. Nocedal & Wright, 1999, §6.2): While some authors report good successes with this method, others note that it is unstable and overly limited. Our Bayesian interpretation adds to the doubts about the SR1 formula, since it identifies it as Gaussian regression with a prior variance involving $V_{i-1}$ with

$$V_{i-1}\boldsymbol{s}_i = (\boldsymbol{y}_i - B_{i-1}\boldsymbol{s}_i), \quad (9)$$

a data-dependent prior covariance. Given the prior (4), there is no rank 1 update rule that gives a symmetric posterior. This blemish of rank-1 updates is also reflected in Eq. (6): Uncertainty drops only in the "row", or "primal" subspace of the belief (the right hand side of the Kronecker product in the covariance). While this still means uncertainty goes toward 0 over time, it does so in an asymmetric way.

---

[1]Here and later, detailed derivations are left out due to space constraints. They can be found in an upcoming journal version of this paper, currently under review.

## 2.1. Symmetric Estimates, but no Symmetric Beliefs

The proper probabilistic way to encode Hessians' symmetry is to include an additional likelihood term

$$\delta(\Delta \vec{B} - \vec{0}) = \lim_{\tau \to 0} \mathcal{N}(\vec{0}, \Delta \vec{B}, \tau \mathbf{I}) \qquad (10)$$

using $\Delta$, the *antisymmetry* operator – the linear map defined through

$$\Delta \vec{X} = \frac{1}{2}\overrightarrow{(X - X^\top)}. \qquad (11)$$

Since this is a linear map, the resulting posterior is analytic, and Gaussian. But the rank of $\Delta$ is $1/2 \cdot N(N-1)$ (e.g. Lütkepohl, 1996, §4.3.1, Eqs. 12 & 20), so the corresponding update rule does not obey the first requirement of Section 1. However, the structure of Eq. (6) hints at another idea, which in fact turns out to give rise to the most popular quasi-Newton methods. We introduce a second, *dual* observation (dual, as in "dual vector space", not as in "primal-dual optimization").

$$\begin{aligned} p(\boldsymbol{y}_i^\top \mid B, \boldsymbol{s}_i^\top) &= \delta(\boldsymbol{y}_i^\top - \boldsymbol{s}_i^\top B) \\ &= \lim_{\gamma \to 0} \mathcal{N}\left[\boldsymbol{y}_i^\top; \boldsymbol{s}_i^\top \vec{B}, \gamma \otimes V_i\right]. \end{aligned} \qquad (12)$$

The posterior after both primal and dual observation is a Gaussian with mean and covariance

$$B_i = B_{i-1} + \frac{(\boldsymbol{y}_i - B_{i-1}\boldsymbol{s}_i)\boldsymbol{s}_i^\top V_{i-1}^\top}{\boldsymbol{s}_i^\top V_{i-1} \boldsymbol{s}_i} + \frac{V_{i-1}\boldsymbol{s}_i(\boldsymbol{y}_i - B_{i-1}\boldsymbol{s}_i)^\top}{\boldsymbol{s}_i^\top V_{i-1}\boldsymbol{s}_i} \\ - \frac{V_{i-1}\boldsymbol{s}_i(\boldsymbol{s}_i^\top(\boldsymbol{y}_i - B_{i-1}\boldsymbol{s}_i))\boldsymbol{s}_i^\top V_{i-1}}{(\boldsymbol{s}_i^\top V_{i-1}\boldsymbol{s}_i)^2} \qquad (13)$$

$$\Sigma_i = \left(V_{i-1} - \frac{V_{i-1}\boldsymbol{s}_i\boldsymbol{s}_i^\top V_{i-1}}{\boldsymbol{s}_i^\top V_{i-1}\boldsymbol{s}_i}\right) \otimes V_i = V_i \otimes V_i. \qquad (14)$$

The posterior mean is clearly symmetric if $B_{i-1}$ is symmetric (as $V_{i-1}$ is symmetric by definition). Choosing the unit prior $\Sigma_{i-1} = \mathbf{I} \otimes \mathbf{I}$ once more, Eq. (13) gives what is known as Powell's (1970) *symmetric Broyden (PSB)* update. Eq. (13) has previously been known to be the most general form of a symmetric rank 2 update obeying the quasi-Newton equation and minimizing a Frobenius regularizer (Dennis & Morée, 1977). This old result is a corollary of our derivations. But note that symmetry only extends to the mean, not the entire belief: In contrast to the posterior generated by Eq. (10), samples from this posterior are, with probability 1, not symmetric. Of course, they can be projected into the space of symmetric matrices by applying the symmetrization operator $\Gamma$ defined by

$$\Gamma \vec{X} = \frac{1}{2}\overrightarrow{(X + X^\top)} \qquad \text{(note that } \mathbf{I} = \Gamma + \Delta; \Gamma\Delta = \mathbf{0}). \qquad (15)$$

Since $\Gamma$ is a symmetric linear operator, the projection of any Gaussian belief $\mathcal{N}(X; X_0, \Sigma)$ onto the space of symmetric matrices is itself a Gaussian $\mathcal{N}(\Gamma X; \Gamma X_0, \Gamma \Sigma \Gamma)$. But symmetrized samples from the posterior of Eqs. (13) & (14) do not necessarily obey the quasi-Newton Equation (2). While Eq. (12) does convey useful information, it is not equivalent to encoding symmetry. It is cheaper, but also weaker, than using the correct likelihood (10).

## 2.2. Positive Definiteness: Meaning or Decoration?

Consider choosing $V_{i-1} = B$. The prior is then

$$p(B) \propto |B|^{-N^2/2} \qquad (16)$$
$$\cdot \exp\left[-\frac{1}{2}\left(N - 2\operatorname{tr}(B_{i-1}B^{-1}) + \operatorname{tr}(B_{i-1}B^{-1}B_{i-1}B^{-1})\right)\right].$$

This is an intriguing prior. Although there is some semblance to the Wishart distribution, the second term in the exponential means this prior is broader than the Wishart. It is not well-defined for degenerate matrices, and it is not clear whether it is proper. It is thus surprising to discover that it engenders the two most popular quasi-Newton methods: If we use the quasi-Newton equation (2) a second time to replace $V_{i-1}\boldsymbol{s} = \boldsymbol{y}$, Eq. (16) gives the *DFP* method (Davidon, 1959; Fletcher & Powell, 1963)

$$B_i = B_{i-1} + \frac{(\boldsymbol{y}_i - B_{i-1}\boldsymbol{s}_i)\boldsymbol{y}_i^\top}{\boldsymbol{s}_i^\top \boldsymbol{y}_i} + \frac{\boldsymbol{y}_i(\boldsymbol{y}_i - B_{i-1}\boldsymbol{s}_i)^\top}{\boldsymbol{y}_i^\top \boldsymbol{s}_i} \\ - \frac{\boldsymbol{y}_i(\boldsymbol{s}_i^\top(\boldsymbol{y}_i - B_{i-1}\boldsymbol{s}_i))\boldsymbol{y}_i^\top}{(\boldsymbol{y}_i^\top \boldsymbol{s}_i)^2}. \qquad (17)$$

And, if we exchange in the entire preceding derivation $\boldsymbol{s} \leftrightarrow \boldsymbol{y}$, $B \leftrightarrow B^{-1}$, $B_{i-1} \leftrightarrow B_{i-1}^{-1}$, then we arrive at the *BFGS* method (Broyden, 1969; Fletcher, 1970; Goldfarb, 1970; Shanno, 1970), which ranks among the most widely used algorithms in machine learning overall. DFP and BFGS owe much of their popularity to the fact that the updated $B_{i,\text{DFP}}$ and $B_{i,\text{BFGS}}^{-1}$ are guaranteed to be positive definite whenever $B_{i-1,\text{DFP}}$ and $B_{i-1,\text{BFGS}}^{-1}$ are positive definite, respectively, and additionally $\boldsymbol{y}_i^\top \boldsymbol{s}_i > 0$. How helpful is this property? It is relatively straightforward to extend a theorem by Dennis & Morée (1977) to find that, assuming $B_{i-1}$ is positive definite, the posterior mean of Eq. (13) is positive definite if, and only if,

$$0 < (\boldsymbol{y}_i^\top B_{i-1}^{-1} V_{i-1} \boldsymbol{s}_i)^2 \qquad (18)$$
$$+ (\boldsymbol{y}_i - B_{i-1}\boldsymbol{s}_i)^\top B_{i-1}^{-1} \boldsymbol{y}_i \cdot \boldsymbol{s}_i^\top V_{i-1} B_{i-1}^{-1} V_{i-1} \boldsymbol{s}_i$$
$$= \boldsymbol{s}_i^\top V_{i-1}[B_{i-1}^{-1}\boldsymbol{y}_i\boldsymbol{y}_i^\top B_{i-1}^{-1} - \boldsymbol{y}^\top B_{i-1}^{-1}\boldsymbol{y}_i + \boldsymbol{s}_i^\top \boldsymbol{y}_i]V_{i-1}\boldsymbol{s}_i.$$

If the prior covariance is not to depend on the data, it is thus impossible to guarantee positive definiteness in this framework – BFGS and DFP circumvent this conceptual issue by choosing $V_{i-1} = B$, then applying Eq. (2) a second time. But, even casting aside such philosophical reservations, our analysis also casts doubt upon the efficacy of the way in which DFP and BFGS achieve positive definiteness: Eq. (16) does not exclude indefinite matrices; in fact it assigns positive measure to every invertible matrix. For example, under a mean $B_{i-1} = \boldsymbol{I}$, the indefinite matrix $B = \text{diag}(1, -1)$ is assigned $p(B) \propto \exp(-2)$. DFP and BFGS achieve positive definiteness, not by including additional information, but by manipulating the prior such that, as if by accident, the *MAP estimator* (not the belief) happens to be positive definite. These observations do not rule out any utility of guaranteeing positive definiteness in this way. But there is less value in the positive definiteness guarantee of DFP and BFGS than previously thought. The algorithm should aim to find the "best" positive definite explanation for the data, not "any" such explanation.

### 2.3. Rank $M$ Updates

The classical quasi-Newton algorithms update the mean of the belief at every step in a rank 2 operation, then, implicitly, reset their uncertainty in the next step, thereby discarding information acquired earlier. Albeit inelegant from a Bayesian point of view, this scheme is still a good idea given other aspects of the framework: Since the quasi-Newton likelihood models the objective function as a quadratic, model mismatch would lead to strong overfitting under exact Bayesian inference. But it is instructive to consider the effect of encoding more than just the most recent observation. It is straightforward to extend Eq. (2) to observations $(Y, S)$ from several line searches:

$$Y_{nm} = \nabla_n f(\boldsymbol{x}_{i_m}) - \nabla_n f(\boldsymbol{x}_{i_m-1})$$
$$S_{nm} = x_{i_m,n} - x_{i_m-1,n} \qquad (19)$$

Given a prior $p(B) = \mathcal{N}(B; B_0, V_0)$, the Gaussian posterior then has mean and covariance

$$B_i = B_0 + (Y - B_0 S)(S^\top V_0 S)^{-1} S^\top V_0 \qquad (20)$$
$$+ V_0 S (S^\top V_0 S)^{-1} (Y - B_0 S)^\top$$
$$- V S (S^\top V_0 S)^{-1} (S^\top (Y - B_0 S))(S^\top V S)^{-1} S^\top V_0$$
$$\Sigma_i = \left(V_0 - V_0 S(S^\top V_0 S)^{-1} S^\top V_0\right) \qquad (21)$$
$$\otimes \left(V_0 - V_0 S(S^\top V_0 S)^{-1} S^\top V_0\right).$$

Here, the absence of information about the symmetry of the Hessian becomes even more obvious: No matter the prior covariance $V_0$, because of the term $S^\top Y$ in the third line of Eq. (20), the posterior mean is not in general symmetric, *unless* $Y = BS$, (e.g. if the objective function is in fact a quadratic).

### 2.4. Summary

The preceding section showed that quasi-Newton algorithms, including the state-of-the-art BFGS and DFP algorithms, can be interpreted as approximate Bayesian regression from the primal and dual likelihood of Eqs. (2) and (12) under varying priors, in the following sense: At each quasi-Newton step, fix a Gaussian prior ad hoc, update the mean, then "forget" the covariance update. Two particularly interesting observations concern the way in which the desiderata of symmetry and positive definiteness of the MAP estimator are achieved in these algorithms. Symmetry is encoded via dual observations, which is a useful but imperfect shortcut. Positive definiteness is achieved not by encoding relevant information, but by shifting the prior post hoc. It is thus doubtable whether the proven good performance of BFGS and DFP is actually down to positive definiteness, instead of a simpler effect of moving from the clearly pathological formulation of Broyden's method to dual observations and a less restrictive (though nontrivial) prior.

## 3. A Nonparametric Bayesian Quasi-Newton Method

Section 2 used the probabilistic perspective to gain novel insight into classical methods. In this second part of the paper we depart from the traditional framework to construct a nonparametric, Bayesian quasi-Newton method, de novo. To motivate this effort, notice some further deficiencies of DFP/BFGS regarding use of available information: Eq. (2) assumes that the function is (locally) a quadratic. Old observations collected "far" from the current location (in the sense that a second order expansion is a poor approximation) may thus be useless or even harmful. The fact that the function is not quadratic should be part of the model. On an only slightly related point, individual line searches typically involve several evaluations of the objective function $f$ and its gradient; but the algorithms only make use of one of those (the last one). This is clearly wasteful, but even the exact Bayesian parametric algorithm of Section 2.3 has this problem: Because a matrix $S$ of several observations along one line search has rank 1, the inverse of $S^\top V_0 S$ is not defined. The following section will address all these issues. Several aspects of the resulting algorithm are involved. Derivations can be found in the journal version. A MATLAB implementation can be found at `www.probabilistic-optimization.org`.

### 3.1. A Nonparametric Prior

Defining a prior for the function $B : \mathbb{R}^N \to \mathbb{R}^{N \times N}$, we choose a set of $N^2$ correlated Gaussian processes. The mean function is assumed to be an arbitrary integrable function $B_0(x)$ (in our implementation we use a constant function, but the analytic derivations do not need to be so restrictive). The core idea is to assume that the covariance between the element $B_{ij}$ at location $x_\triangledown$ and the entry $B_{k\ell}$ at location $x_\triangle$ is

$$\operatorname{cov}(B_{ij}(x_\triangledown), B_{k\ell}(x_\triangle)) = k_{ik}(x_\triangledown^\top, x_\triangle^\top) k_{j\ell}(x_\triangledown, x_\triangle) \\ = (\boldsymbol{k} \otimes \boldsymbol{k})_{(ij)(k\ell)}(x_\triangledown, x_\triangle) \quad (22)$$

with an $N \times N$ matrix of kernels, $\boldsymbol{k}$. To give a more concrete intuition: In our implementation we use one joint squared exponential kernel for all elements. I.e.

$$k_{ij}(x_\triangledown, x_\triangle) = V_{ij} \exp\left(-\frac{1}{2}(x_\triangledown - x_\triangle)^\top \Lambda^{-1}(x_\triangledown, x_\triangle)\right) \quad (23)$$

with a positive definite matrix $V$ and length scales $\Lambda$. Other kernels can of course be chosen; but it will become clear that an important practical requirement is the ability to efficiently integrate the kernel. This is feasible, though nontrivial, with the squared exponential kernel. Another option, not yet explored by us, may be offered by spline kernels (Minka, 2000).

### 3.2. Line Integral Observations

For the Hessian $B(x)$ of a general function $f$, the quasi-Newton equation (2) is only a zeroth order approximation (a second-order approximation to $f$ itself), assuming a constant Hessian everywhere. In our treatment, we will replace it with the exact statement: We observe the value of the *line integral* along the path $r^i : [0,1] \to \mathbb{R}^N$, $r^i(t) = x_{i-1} + t(x_i - x_{i-1})$.

$$Y_{ni} = \sum_m \int_{r^i_m} B_{nm}(x)\,dx_m = \sum_m S_{mi} \int_0^1 B_{nm}(r^i(t))\,dt. \quad (24)$$

This uses the classic result that line integrals over scalar fields, such as $B(x)$, are fully defined by the path's start and end point, irrespective of the path itself. Hence, the nonparametric version of the quasi-Newton equation is the likelihood

$$p(Y \mid B(x), \mathfrak{S}) = \lim_{\beta \to 0} \mathcal{N}\left[Y; \mathfrak{S}^\top \overrightarrow{B}, k \otimes \beta \boldsymbol{I}_M\right] \quad (25)$$

with a linear operator ($\odot$ denotes the Hadamard, or element-wise product $(a \odot b)_{k\ell} = a_{k\ell} b_{k\ell}$)

$$\mathfrak{S} = \boldsymbol{I} \otimes \left(\int_0^1 dt \odot S\right). \quad (26)$$

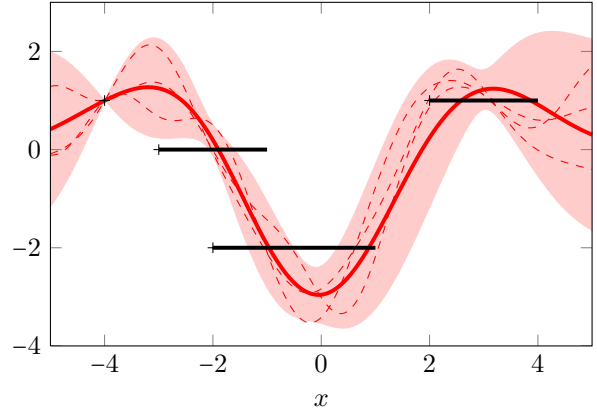

*Figure 1.* One-dimensional Gaussian process inference from integral observations (squared exponential kernel). Four observations, average values (integral value divided by length of integration region) and integration regions denoted by black bars. Posterior mean in thick red, two standard deviations as shaded region, three samples as dashed lines. The left-most integral is over a very small region, so that it essentially reduces to the classical case of a local observation. Corresponding integrals over the mean, and each sample, are consistent with the integral observations.

### 3.3. Gaussian Process Inference from Integral Observations

Because the Gaussian exponential family is closed under linear transformations, Gaussian process inference is analytic under any linear operator. Since integration is a linear operation, Gaussian process inference is possible, in closed form, from integral observations. Nevertheless, this idea has only rarely been used in the literature (e.g. by Minka, 2000). Figure 1 gives a 1D toy example for intuition.

The posterior distribution under our nonparametric prior, the likelihood of Eq. (25) and its dual equivalent is a Gaussian process with mean and covariance functions

$$B_\diamond(x_\triangledown) = B_0(x_\triangledown) + (Y - \mathfrak{B}_0)\mathfrak{K}^{-1}\mathfrak{k}^\top(x_\triangledown) \\ + \mathfrak{k}(x_\triangledown)\mathfrak{K}^{-1}(Y - \mathfrak{B}_0)^\top \quad (27) \\ - \mathfrak{k}(x_\triangledown)\mathfrak{K}^{-1}S^\top(Y - \mathfrak{B}_0)\mathfrak{K}^{-1}\mathfrak{k}^\top(x_\triangledown)$$

$$\Sigma_\diamond(x_\triangledown, x_\triangle) = \left[k(x_\triangledown^\top, x_\triangle^\top) - \mathfrak{k}(x_\triangledown^\top)\mathfrak{K}^{-1}\mathfrak{k}^\top(x_\triangle)\right] \quad (28) \\ \otimes \left[k(x_\triangledown, x_\triangle) - \mathfrak{k}(x_\triangledown)\mathfrak{K}^{-1}\mathfrak{k}^\top(x_\triangle)\right].$$

This uses $\mathfrak{B}_0 \in \mathbb{R}^{N \times M}$, the function $\mathfrak{k} : \mathbb{R} \to \mathbb{R}^{N \times M}$ and the Gram matrix $\mathfrak{K} \in \mathbb{R}^{M \times M}$, defined by

$$\mathfrak{B}_{0,nm} = \sum_{\ell} S_{\ell m} \int_0^1 B_{0,n\ell}(r^\ell(t))\,\mathrm{d}t$$

$$\mathfrak{k}_{nm}(x_\triangledown) = \sum_{\ell} S_{\ell m} \int_0^1 k(x_\triangledown, r^\ell(t))\,\mathrm{d}t \qquad (29)$$

$$\mathfrak{K}_{pq} = \sum_{\ell,j} S_{\ell p} S_{jq} \iint_0^1 k(r^\ell(t), r^j(t'))\,\mathrm{d}t\,\mathrm{d}t'.$$

These objects are homologous to concepts in canonical Gaussian process inference: $\mathfrak{B}_{0,nm}$ is the $n$-th mean prediction along the $m$-th line integral observation. $\mathfrak{k}_{nm}(x_\triangledown)$ is the covariance between the $n$-th column of the Hessian at location $x_\triangledown$ and the $m$-th line-integral observation. $\mathfrak{K}_{pq}$ is the covariance between the $p$-th and $q$-th line integral observations. An important aspect is that, because $k$ is a positive definite kernel, unless two observations are exactly identical, $\mathfrak{K}$ has full rank $M$ (the number of function evaluations), even if several observations take place within one shared 1-dimensional subspace. So it is possible to make full use of *all* function evaluations made during line searches, not just the first and last one, as in the classical setting. A downside is that evaluating the mean function involves finding the inverse of $\mathfrak{K}$, at cost $\mathcal{O}(M^3)$. Two aspects of numerical optimization make this issue less problematic than one might think. First, solving an optimization problem takes finite time, often just a few hundred evaluations; so the cubic cost in $M$ is often manageable. Where it is not, note that, because optimization proceeds along a trajectory through the parameter space, old observations tend to have low covariance with the Hessian at the current location, and thus a small effect on the local mean estimate (the effect of this influence is measured by $\mathfrak{k}\mathfrak{K}^{-1}$). So they can often simply be ignored.

### 3.4. Numerical Implementation

As mentioned above, for a concrete implementation, we chose to use the squared exponential kernel (23), and a constant mean function assigning $B_0(x_\triangledown) = \boldsymbol{I}$ everywhere. It is another advantage of the Bayesian formulation that prior assumptions are easy to analyze and understand: The squared exponential prior amounts to the assumption that the elements of the Hessian vary independently over the parameter space, on one unique set of length-scales $\Lambda$. Multiple length scales could be modeled using sums of kernels, but our implementation does not currently offer this option.

Changing the length scales $\Lambda$ amounts to *automatic pre-conditioning*, another benefit of a Bayesian formulation that we cannot dwell on for space reasons. Hyperparameters could be fitted by type-II maximum likelihood, as in canonical Gaussian process regression. Unfortunately, this is an optimization problem itself. Another option is to instead fix the hyperparameters ad hoc by tracking the signal variance to fix $V$ in Eq. (23) and the relative change along line searches to fix $\Lambda$.

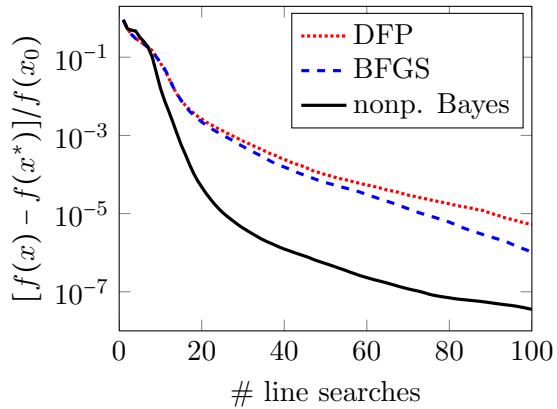

*Figure 2.* Minimizing the logarithm of a 200-dimensional product of Gamma distributions. Averages over 20 sampled problems; plotted is the relative distance from initial function value (shared by all algorithms) to the minimum, as a function of the number of line searches (all algorithms use the same line search method).

Implementing the integrals of Eq. (29) for the squared-exponential kernel, particularly those in $\mathfrak{K}$, is nontrivial, because definite integrals over Gaussians are not analytic. $\mathfrak{k}$ involves the error function, for which good double-precision approximations are widely available. The integrals in $\mathfrak{K}$ are of two distinct types: The covariance between observations made as part of the same line search involve 1D integrals of the error function, which can be analytically reduced to the error and exponential functions[2]. The covariance between observations made during different line searches are bivariate Gaussian integrals. Fortunately, good, light-weight numerical approximations are available for this problem (Genz, 2004).

From Sec. 1, recall that updating the search direction requires the *inverse* of $B$. Explicit inversion costs $\mathcal{O}(N^3)$, but the inverse can be constructed analytically, from the matrix inversion lemma, in $\mathcal{O}(N^2 + NM + M^3)$. Using a diagonal prior mean $B_0$ and an argument largely analogous to the derivation of the L-BFGS algorithm (Nocedal, 1980) lowers cost to $\mathcal{O}(NM + M^3)$, linear in $N$. The nonparametric method is thus applicable to problems of even very high dimensionality.

---
[2]Jaakko Peltonen, 2011, personal communication

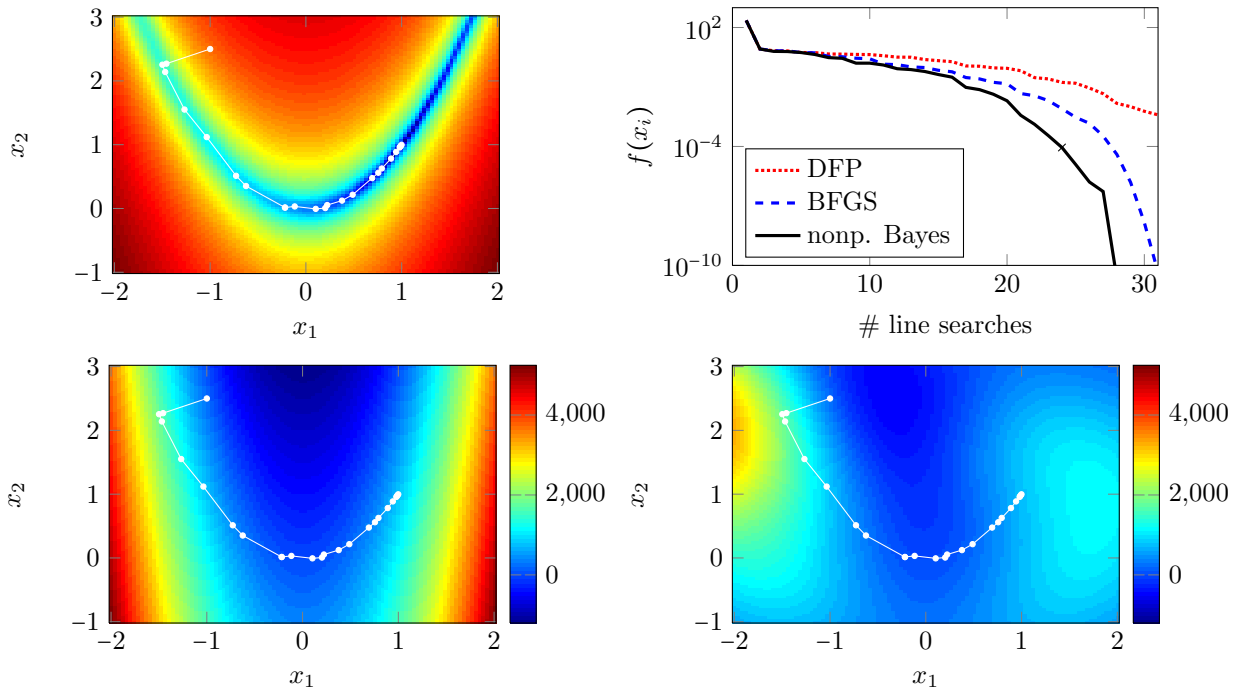

*Figure 3.* Minimizing Rosenbrock's polynomial, a non-convex function which has its only minimum at (1,1). All algorithms start from (-1,2.5). **Top left:** Function values, line search trajectory of the Bayesian algorithm in white. **Top right:** function value as function of the number of line searches. The black cross marks the point where the Bayesian method switches to a local parametric model for numerical stability. **Bottom Left:** True value of the (1,1) element of the Hessian (other elements have less interesting structure). **Bottom Right:** Mean estimate of the Bayesian regressor, showing good agreement in the regions visited by the algorithm (corresponding uncertainty measure not shown).

## 4. Experiments

Figure 2 shows averages of experiments on a 200-dimensional domain. The objective functions were the logarithms of products of Gamma distributions with different parameters for each dimension (a simplified version of hyperparameter learning for Gaussian process regression). In this experiment, the nonparametric algorithm outperforms its predecessors strongly. The performance advantage is not always so drastic (the journal version contains additional empirical results, including less pronounced cases). Despite the relatively precise numerical treatment of the integrals involved, the nonparametric Bayesian quasi-Newton algorithm poses more numerical challenges than its predecessors. This issue becomes clear when minimizing quadratic functions, whose constant Hessian voids the modelling advantage of the nonparametric method: The Bayesian algorithm behaves more regularly initially, but towards the end of the optimization process the numerical conditioning of the Bayesian algorithms begins to play a role, offering an advantage to the better conditioned older methods. At this small scale, however, the Hessian is essentially constant, and the function is well described by a local model. In our practical implementation, we check for convergence, then pass the learned inverse Hessian to the better conditioned BFGS for the final few steps.

An additional benefit of the nonparametric formulation is the availability of a global estimate of the Hessian function. Figure 3 illustrates this point with results from a popular two-dimensional test problem – Rosenbrock's polynomial (details in caption). This figure is mostly for intuition: Rosenbrock's valley is challenging even for the exact Newton method since it breaks the line search paradigm, so the similarity between the methods on this problem is not particularly indicative of general performance.

**Cost**  As pointed out above, the computational complexity of this algorithm, given a diagonal prior mean, is $\mathcal{O}(NM + M^3)$ per update of the search direction, where $M$ is the number of function evaluations used to build the model (which can be controlled ad hoc within the algorithm by excluding redundant or irrelevant evaluations). This compares to $\mathcal{O}(NM)$ for the corresponding cases of DFP and BFGS. Although the overhead created by the squared-exponential integrals is nontrivial, we found the computational demands of our implementation manageable: In our experiments, the cost of constructing and inverting the matrix $\mathfrak{K}$ was negligible, and could, in very time-sensitive settings, be further reduced by a more efficient implementation.

## 5. Outlook

Owing to the limitations of a conference publication, we have only outlined many of our core results. To give an intuition for the potential of probabilistic formulations of numerical optimization, consider some of the most immediate future work: Perhaps the most obvious insight is that Gaussian process integration is trivial to extend to noisy evaluations. In combination with a robust replacement for the traditional line searches, our work may thus lead to robust numerical optimizers. Repeated integration, and non-Gaussian likelihoods in combination with approximate inference, may allow optimization without gradients, and from only gradient sign observations, respectively. Structured and hierarchical priors are a third direction, offering new avenues for optimization of very high-dimensional functions.

## 6. Conclusion

We have shown that the most popular quasi-Newton algorithms can be interpreted as approximations to Bayesian regression under Gaussian and other priors. This deepens our understanding of these algorithms. In particular, it emerged that symmetry in the estimators of SR1, PSB, DFP and BFGS, and positive definiteness in those of DFP and BFGS, are encoded in only approximate, incomplete ways.

As a parallel result, our analysis also gives rise to a new class of Bayesian nonparametric quasi-Newton algorithms. These use a kernel model to utilize all observations in each line-search, explicitly track uncertainty, and thus achieve faster convergence towards the true Hessian. While the new methods are not trivial to understand and implement, their computational cost lies within a constant of that of their predecessors. A demonstrative implementation can be found at `www.probabilistic-optimization.org`.

## Acknowledgments

The authors thank Christian Schuler, Tom Minka and Carl Rasmussen for helpful discussions, as well as Carl Rasmussen for his release of `minimize.m`, which simplified development. MK is supported by a grant from Microsoft Research Ltd.